\documentclass[final]{article}
\usepackage{graphicx}
\usepackage{amsmath}
\usepackage{amsfonts}
\usepackage{amssymb}

\usepackage{graphicx}
\usepackage{amsmath}
\usepackage{amsfonts}
\usepackage{amssymb}

\newcommand{\interval}{I}
\renewcommand{\deg}{n} % Polynomial degree
\newcommand{\cont}{m}  % Continuity order
\newcommand{\dif}{\mathrm d}
\newcommand{\dx}{\dif x}
\newcommand{\dt}{\dif t}
\let\phi\varphi

\renewcommand{\P}{\mathbb P}
\newcommand{\R}{\mathbb R}
\newcommand{\nodal}{\mathcal N}
\newcommand{\cint}{I}
\newcommand{\wint}{\Pi}
\newcommand{\ordertens}{N}

\newtheorem{theorem}{Theorem}
\newtheorem{lemma}[theorem]{Lemma}
\newtheorem{remark}[theorem]{Remark}

\title{A Tensor-Product Finite Element Cochain Complex with Arbitrary Continuity}

\author{Francesca Bonizzoni\\%\thanks{%
  Department of Mathematics\\
  University of Augsburg\\
\texttt{francesca.bonizzoni@math.uni-augsburg.de}    
\and Guido Kanschat\\%\thanks{%
Interdisciplinary Center for Scientific Computing (IWR)\\
Heidelberg University\\
\texttt{kanschat@uni-heidelberg.de}
}

\begin{document}
\maketitle
\abstract{
We develop tensor product finite element cochain complexes of arbitrary smoothness on Cartesian meshes of arbitrary dimension. The first step is the construction of a one-dimensional $C^m$-conforming finite element cochain complex based on a modified Hermite interpolation operator, which is proved to commute with the exterior derivative by means of a general commutation lemma. Adhering to a strict tensor product construction we then derive finite element complexes in higher dimensions.
}

%%%%%%%%%%%%%%%%%%%%%%%%%%%%%%%%%%%%%%%%%%%%%%%%%%%%%%%%%%%%%%%%%%%%%%%%%%%%%%%%%%%%%%%%%%%%%%%%%%%%%%%%%%%
\section{Introduction}
%%%%%%%%%%%%%%%%%%%%%%%%%%%%%%%%%%%%%%%%%%%%%%%%%%%%%%%%%%%%%%%%%%%%%%%%%%%%%%%%%%%%%%%%%%%%%%%%%%%%%%%%%%%

This article constructs a family of finite element cochain complexes (FECC) of differential forms in any dimension with arbitrary continuity $m$ on Cartesian meshes. The building block is a one-dimensional FECC based on modified Hermite interpolation. It is $C^m$-conforming, meaning that, for all $u\in \mathbb P_n\Lambda^0\subset C^m\Lambda^0$, there holds $\dif u\in\mathbb P_{n-1}\Lambda^1\subset C^{m-1}\Lambda^1$. The FECC in arbitrary dimension is then derived by adhering to a strict tensor product construction.

Finite element cochain complexes have developed as an important tool to construct and analyze discretization methods for vector fields and operators of the de Rham complex, see e.~g.~\cite{ArnoldFalkWinther}. By applying the Bernstein-Gelfand-Gelfand construction, they can also be used to obtain mixed discretization schemes for more complex situations like elasticity~\cite{ArnoldHu}. There, a sequence of de Rham sequences of varying continuity is employed.

A first step in the direction of FECC with arbitrary continuity was accomplished in~\cite{BonizzoniKanschat21}, where the authors introduce finite element (FE) spaces forming exact $H^1$-conforming cochain subcomplexes of the de Rham complex. The present contribution generalizes~\cite{BonizzoniKanschat21} from $H^1$-conformity to arbitrary regularity.

A fundamental property of an FECC is the commutation between the exterior derivative and the FE interpolation operator, since it immediately implies that the FECC is a subcomplex of the continuous complex. We prove this commutativity using a general lemma, which we introduce in Section 2. In Section 3, we develop FEs that admit this commutation property in one dimension. Tensor products of these elements are used in Section 4 to construct the FECC in higher dimensions.

%%%%%%%%%%%%%%%%%%%%%%%%%%%%%%%%%%%%%%%%%%%%%%%%%%%%%%%%%%%%%%%%%%%%%%%%%%%%%%%%%%%%%%%%%%%%%%%%%%%%%%%%%%%
\section{Cochain Complexes of Finite Element Differential Forms}
\label{sec:cochain_complexes_of_FE_differential_forms}
%%%%%%%%%%%%%%%%%%%%%%%%%%%%%%%%%%%%%%%%%%%%%%%%%%%%%%%%%%%%%%%%%%%%%%%%%%%%%%%%%%%%%%%%%%%%%%%%%%%%%%%%%%%

The construction of FE spaces of differential forms follows the general construction principle of FE spaces, namely first subdividing the domain of computation into a mesh of simple mesh cells ---typically simplices or hypercubes--- and then defining a shape function space on each cell. The continuity between mesh cells is established by the choice of node functionals at the interfaces between cells. These node functionals must also be chosen so that each function in the shape function space is uniquely determined by the values of the node functionals, the so called unisolvence.

The construction of an FECC differs only by the fact that the shape function spaces and node functionals are not chosen independently, but in relation to each other. Therefore, let
\begin{gather*}
    0 \xrightarrow{\;\;\subset \;\;}
    \cdots
    \xrightarrow{\;\; \dif \;\;}
    \P\Lambda^k
    \xrightarrow{\;\; \dif \;\;}
    \P\Lambda^{k+1}
    \xrightarrow{\;\; \dif \;\;}
    \cdots
    \xrightarrow{\;\;0 \;\;}0
\end{gather*}
be an exact sequence of differential forms with polynomial coefficients forming the shape function spaces.
Let $\dif\colon\P\Lambda^k\to\P\Lambda^{k+1}$ denote the exterior derivative, and define $r=\operatorname{rank} \dif$ as well as $n_\nu = \dim \P\Lambda^\nu$ for $\nu=k,k+1$. Now assume that on $\P\Lambda^k$ sets of basis functions $\{\phi^k_j\}_{j=1,\dots,n_k}$
and node functionals $\{\nodal^k_i\}_{i=1,\dots,n_k}$ are chosen with the following properties:
\begin{subequations}
\label{eq:separation-domain}
\begin{xalignat}3
    \label{eq:separation-domain-1}
    \phi^k_j &\in \ker \dif & && j&= r+1,\dots,n_k,\\
    \nodal^k_i(\phi) &= 0 & \phi&\in\ker\dif & i&=1,\dots,r,\\
    \nodal^k_i(\phi) &= 0 &\phi&\in\operatorname{span}\{\phi^k_1,\dots,\phi^k_r\} &i&=r+1,\dots,n_k.
\end{xalignat}
\end{subequations}
Namely, the last $n_k-r$ basis functions span the kernel of $\dif$ and the node functionals separate between kernel and cokernel in the sense that the matrix of node functionals applied to basis functions has the structure
\begin{gather}
    \label{eq:block-mk}
    M_k := \bigl(\nodal_i^k(\phi_j^k)\bigr)_{i,j=1,\dots,n_k}
    =\left(\begin{array}{c|c}
         A_k & 0 \\
         \hline
         0 & B_k
    \end{array}\right).
\end{gather}
For the space $\P\Lambda^{k+1}$ we assume that sets of basis functions $\{\phi^{k+1}_j\}_{j=1,\dots,n_{k+1}}$
and node functionals $\{\nodal^{k+1}_i\}_{i=1,\dots,n_{k+1}}$ are chosen so that $\phi^{k+1}_1\dots,\phi^{k+1}_r$ span the range of $\dif$, and the node functionals separate the range of $\dif$ and its complement, namely
\begin{subequations}
\label{eq:separation-range}
\begin{xalignat}3
    \operatorname{range}\dif & = \operatorname{span}\{\phi^{k+1}_1,\dots,\phi^{k+1}_r\},&&&&\\
    \label{eq:separation-range-1}
    \nodal^{k+1}_i(\phi) &= 0 &\phi&\in\operatorname{range}\dif
    &i&=r+1,\dots,n_{k+1},\\
    \label{eq:separation-range-2}
    \nodal^{k+1}_i(\phi) &= 0 &\phi&\in\operatorname{span}\{\phi^{k+1}_{r+1},\dots,\phi^{k+1}_{n_{k+1}}\} &i&=1,\dots,r.
\end{xalignat}
\end{subequations}

Thus, we have the corresponding structure
\begin{gather*}
    M_{k+1} := \bigl(\nodal_i^{k+1}(\phi_j^{k+1})\bigr)_{i,j=1,\dots,n_{k+1}}
    =\left(\begin{array}{c|c}
         A_{k+1} & 0 \\
         \hline
         0 & B_{k+1}
    \end{array}\right).
\end{gather*}
Note though that the subspaces of the two partitionings are defined in a different way. Nevertheless, both choices are always possible by basic results of linear algebra.

\begin{lemma}
\label{lem:commuting-lemma}
Assume that the shape function bases and node functionals for $\P\Lambda^k$ and $\P\Lambda^{k+1}$ are chosen according to conditions~\eqref{eq:separation-domain} and~\eqref{eq:separation-range}. Assume further that
  \begin{gather}
    \label{eq:commute-basis}
      \begin{aligned}
        \dif\phi^k_j &= \phi^{k+1}_j & \qquad j&=1,\dots,r,\\
        \nodal^{k+1}_i(\dif u) &= \nodal^k_i(u)& \qquad i&=1,\dots,r,
      \end{aligned}
  \end{gather}
  for all $u\in C^\infty\Lambda^k$.
  Then, the FE interpolation operators commute with the exterior derivative $\dif$, that is,
  \begin{gather}
    \label{eq:interp}
    \dif \cint_k u = \cint_{k+1} \dif u \qquad \forall u \in C^\infty\Lambda^k.
  \end{gather}
\end{lemma}

\emph{Proof: } Given the node functionals, the FE interpolation operators $\cint_\nu$, for $\nu=k,k+1$, are defined as
\begin{equation}
    \label{eq:cint}
    \cint_\nu u =\sum_{i=1}^{n_\nu} \nodal^\nu_i(u)
    \left(\sum_{j=1}^{n_\nu} \alpha^\nu_{ij}\phi^\nu_j\right)
    \qquad \forall u\in  C^\infty\Lambda^\nu,
\end{equation}
where the coefficients $\alpha^\nu_{ij}$ are the entries of the matrix $M_\nu^{-1}$. We thus have by~\eqref{eq:separation-range-1}
\begin{gather*}
    \cint_{k+1}(\dif u) = \sum_{i=1}^{r} \nodal^{k+1}_i(\dif u)
    \left(\sum_{j=1}^{r} \alpha^{k+1}_{ij}\phi^{k+1}_j\right)
\end{gather*}
where we shortened the second sum to $r$ terms using the block structure of $M_{k+1}$.
Similarly, by~\eqref{eq:separation-domain-1}
\begin{gather*}
    \dif \cint_k u
    = \sum_{i=1}^{r} \nodal^{k}_i(u)
    \left(\sum_{j=1}^{r} \alpha^{k}_{ij} \dif\phi^{k}_j\right).
\end{gather*}
Thus, applying the commutation conditions in~\eqref{eq:commute-basis}, the commutation result~\eqref{eq:interp} will hold if we can show that there holds
\begin{gather*}
    \alpha_{ij}^k = \alpha_{ij}^{k+1},\qquad i,j=1,\dots,r,
\end{gather*}
namely, if $A_k = A_{k+1}$. Applying again~\eqref{eq:commute-basis},
we have
\begin{gather*}
   A_{k+1}
   = \bigl(\nodal_i^{k+1}(\phi_j^{k+1})\bigr)_{i,j=1,\dots,r}
   = \bigl(\nodal_i^{k+1}(\dif \phi_j^{k})\bigr)_{i,j=1,\dots,r}
   = \bigl(\nodal_i^{k}(\phi_j^{k})\bigr)_{i,j=1,\dots,r}
   = A_{k},
\end{gather*}
hence the commutation property~\eqref{eq:interp} is proved.
\hfill $\square$

%%%%%%%%%%%%%%%%%%%%%%%%%%%%%%%%%%%%%%%%%%%%%%%%%%%%%%%%%%%%%%%%%%%%%%%%%%%%%%%%%%%%%%%%%%%%%%%%%%%%%%%%%%%
\section{The Finite Element Cochain Complex in One Dimension}
%%%%%%%%%%%%%%%%%%%%%%%%%%%%%%%%%%%%%%%%%%%%%%%%%%%%%%%%%%%%%%%%%%%%%%%%%%%%%%%%%%%%%%%%%%%%%%%%%%%%%%%%%%%

A polynomial cochain complex in one dimension has the form
\begin{gather}
    \label{eq:1D_complex}
    0 \xrightarrow{\;\;\subset \;\;}\R
    \xrightarrow{\;\; \subset \;\;}
    \P_{\deg}\Lambda^0
    \xrightarrow{\;\; \dif \;\;}
    \P_{\deg-1}\Lambda^1
    \xrightarrow{\;\; 0 \;\;}
    0.
\end{gather}
Introducing vector proxies, a natural correspondence can be established between the spaces of polynomial 0-forms $\P_{\deg}\Lambda^0$ and polynomial functions $\P_{\deg}$ of degree $\deg$, as well as the spaces of polynomial 1-forms $\P_{\deg-1}\Lambda^1$ and polynomial functions $\P_{\deg-1}$ of degree $\deg -1$. In this view, the exterior derivative $\dif$ mapping $0$-forms to $1$-forms corresponds to the derivative in one dimension. The spaces fulfill the property $\dif\P_{\deg}\Lambda^0\subset\P_{\deg-1}\Lambda^1$.
In order to construct a FE method of such forms, we must complement the spaces of polynomial 0- and 1-forms by suitable node functionals which serve to establish continuity of piecewise polynomials across the interfaces between intervals. We note that a version of~\eqref{eq:1D_complex} for $C^m$-regular functions reads
\begin{gather}
    \label{eq:1D_cont_complex}
    0 \xrightarrow{\;\;\subset \;\;}\R
    \xrightarrow{\;\; \subset \;\;}
    C^{\cont}\Lambda^0
    \xrightarrow{\;\; \dif \;\;}
    C^{\cont-1}\Lambda^1
    \xrightarrow{\;\; 0 \;\;}
    0.
\end{gather}
Thus, piecewise polynomials in $C^\cont$ for 0-forms must be paired with such in $C^{\cont-1}$ in order to be consistent. In the next sections, we establish such continuity for polynomial differential forms on the interval $\interval=[0,1]$ by a mixture of node functionals for Hermite interpolation and integral moments.

%%%%%%%%%%%%%%%%%%%%%%%%%%%%%%%%%%%%%%%%%%%%%%%%%%%%%%%%%%%%%%%%%%%%%%%%%%%%%%%%%%%%%%%%%%%%%%%%%%%%%%%%%%%
\subsection{Finite elements for 0-forms}
\label{sec:fem0}
%%%%%%%%%%%%%%%%%%%%%%%%%%%%%%%%%%%%%%%%%%%%%%%%%%%%%%%%%%%%%%%%%%%%%%%%%%%%%%%%%%%%%%%%%%%%%%%%%%%%%%%%%%%

We equip the space $\P_{\deg}(\interval)$ with the following set of node functionals $\{\nodal^0_i\}_{i=1,\ldots,\deg+1}$. First, we use the standard derivative degrees of freedom of two-point Hermitian interpolation:
\begin{gather}
    \label{eq:nodes-1}
 \begin{array}{rl}
    \nodal^0_{2i+1} (u) &= \partial_x^{i+1} u(0)\\
    \nodal^0_{2i+2} (u) &= \partial_x^{i+1} u(1)
\end{array}
    \quad i=0,\ldots,\cont-1.
\end{gather}
They are complemented by interior moments
\begin{equation}
    \label{eq:nodes-4}
    \nodal^0_{2m+i} (u) = \int_0^1 \ell_{i-1} \partial_x u\,\dx,
    \qquad i=1,\dots,\deg-2m,
\end{equation}
where $\ell_i\in\P_i(\interval)$ denotes the Legendre polynomial of degree $i$
on the interval $[0,1]$,
normalized with the condition $\ell_i(1)=1$. Finally, the last degree of freedom is
\begin{align}
    \label{eq:nodes-3}
    \nodal^0_{\deg+1}(u) &= u(1)+u(0).
\end{align}
Note that in~\eqref{eq:nodes-4}
\begin{align}
    \label{eq:nodes-2}
    \nodal^0_{2m+1}(u) &=  \int_0^1 \partial_x u\,\dx
    =u(1)-u(0).
\end{align}
Hence, $\nodal^0_{2m+1}$ and $\nodal^0_{n+1}$ are linear combinations of the function value degrees of freedom of Hermite interpolation, which together with~\eqref{eq:nodes-3} ensure continuity at the end points of the interval $\interval$.

For $\alpha=0,1$ and $\beta=0,\ldots,m$ let
$h_{\alpha,\beta}(x)$ 
be the standard polynomial basis used for Hermite interpolation such that for $\beta,\gamma=0,\dots,m$
there holds
\begin{xalignat*}2
\partial_x^\gamma h_{0,\beta}(0) &= \delta_{\beta\gamma},&
\partial_x^\gamma h_{0,\beta}(1) &= 0,\\
\partial_x^\gamma h_{1,\beta}(1) &= \delta_{\beta\gamma},&
\partial_x^\gamma h_{1,\beta}(0) &= 0.
\end{xalignat*}
We obtain basis functions for the space of polynomial 0-forms as follows: first,
\begin{xalignat*}3
    \phi^0_{2j+1} &= h_{0,j+1},
    &\phi^0_{2j+2} &= h_{1,j+1},
    & j&=0,\dots,m-1,\\
    \phi^0_{2m+1} &= \tfrac12(h_{1,0}-h_{0,0}),
    &
    \phi^0_{n+1} &= \tfrac12(h_{0,0}+h_{1,0}) \equiv \tfrac12.
\end{xalignat*}
We refer to this set of basis functions as the subset of Hermite interpolating basis functions. Except for $\phi^0_{n+1}$, they are all polynomials of degree $2m+1$. Note in particular that $\phi^0_{2m+1}$ has been defined following the alternative expression of $\nodal^0_{2m+1}$ given in~\eqref{eq:nodes-2}.

Let $L^\alpha_j(x)$ denote the polynomials obtained by subsequent integration of Legendre polynomials. They are recursively defined by $L^0_j = \ell_j$ and
\begin{gather*}
    L^{\alpha+1}_j(x) = \int_0^x L^\alpha_j(t)\dt.
\end{gather*}
The remaining basis functions are then defined as
\begin{gather}
    \label{eq:basis-descending}
    \phi^0_{2m+j} = L^{m+1}_{j+m-1}, \qquad j=2,\ldots,n-2m.
\end{gather}
Note that $\phi^0_{2m+j}$ has an $m$-fold root at both interval ends and is of polynomial degree $2m+j$.

%%%%%%%%%%%%%%%%%%%%%%%%%%%%%%%%%%%%%%%%%%%%%%%%%%%%%%%%%%%%%%%%%%%%%%%%%%%%%%%%%%%%%%%%%%%%%%%%%%%%%%%%%%%
\subsection{Finite elements for 1-forms}
\label{sec:fem1}
%%%%%%%%%%%%%%%%%%%%%%%%%%%%%%%%%%%%%%%%%%%%%%%%%%%%%%%%%%%%%%%%%%%%%%%%%%%%%%%%%%%%%%%%%%%%%%%%%%%%%%%%%%%

We construct node functionals and basis functions for 1-forms on $\interval$ such that the commuting interpolation result of the Section~\ref{sec:cochain_complexes_of_FE_differential_forms} holds. First, we observe that according to~\eqref{eq:1D_complex} the dimension is reduced by one. Thus, we define basis functions for $\P_{n-1}\Lambda^1$ by
\begin{gather*}
    \phi^1_j = \dif \phi^0_j,\qquad j=1,\dots,n.
\end{gather*}
Note that $\phi^0_{n+1}$ is the constant function which is mapped to zero by the derivative. The linear independence of the basis functions will follow from the unisolvence proven in the next section.

Similarly, we define node functionals with commutation in mind: first, let
\begin{gather}
    \label{eq:nodes-5}
 \begin{array}{rl}
    \nodal^1_{2i+1} (v) &= \partial_x^i v(0)\\
    \nodal^1_{2i+2} (v) &= \partial_x^i v(1)
\end{array}
    \quad i=0,\ldots,\cont-1.
\end{gather}
Hence, we have for any $u\in\P_n\Lambda^0$:
\begin{gather*}
    \nodal^1_{i} (\dif u) = \nodal^0_{i} (u)
    \quad i=1,\ldots,2\cont.
\end{gather*}
The remaining node functionals are
\begin{gather}
        \label{eq:nodes-6}
    \nodal^1_{2m+i} (v) = \int_0^1 \ell_{i-1} v\, \dx,
    \qquad i=1,\dots,\deg-2m.
\end{gather}
Comparing with~\eqref{eq:nodes-4}, we obtain for any $u\in\P_n\Lambda^0$:
\begin{gather*}
    \nodal^1_{i} (\dif u) = \nodal^0_{i} (u)
    \quad i=2m+1,\ldots,n.
\end{gather*}
A particular consequence of this construction is
\begin{gather}
    \label{eq:nodal-matrix-commute}
    \nodal^1_i(\phi^1_j) = \nodal^0_i(\phi^0_j),
    \qquad i,j=1,\dots,n.
\end{gather}

%%%%%%%%%%%%%%%%%%%%%%%%%%%%%%%%%%%%%%%%%%%%%%%%%%%%%%%%%%%%%%%%%%%%%%%%%%%%%%%%%%%%%%%%%%%%%%%%%%%%%%%%%%%
\subsection{Unisolvence}
%%%%%%%%%%%%%%%%%%%%%%%%%%%%%%%%%%%%%%%%%%%%%%%%%%%%%%%%%%%%%%%%%%%%%%%%%%%%%%%%%%%%%%%%%%%%%%%%%%%%%%%%%%%

In order to prove unisolvence of the FE for $\P_n\Lambda^0$, we study the structure of nonzero entries of the matrix $a_{ij} = \nodal^0_i(\phi^0_j)$. To this end, we split node functionals and basis functions into three groups:
\begin{enumerate}
    \item Hermite interpolation: basis functions and node functionals with index $i=1,\ldots,2m+1$.
    \item Additional basis functions (bubble functions) and node functionals with indices $2m+2$ to $n$.
    \item The basis function in the kernel of $\dif$ and according node functional with index $n+1$.
\end{enumerate}
By construction, there holds immediately
\begin{xalignat*}2
     \nodal^0_i(\phi^0_j) &= \delta_{ij}, & i,j&=1,\dots, 2m+1,\\
     \nodal^0_i(\phi^0_{n+1}) &= \delta_{i,n+1}, & i&=1,\dots,n+1,\\
     \nodal^0_{n+1}(\phi^0_j) &= \delta_{j,n+1}, &
     j&=1,\dots,n+1.
\end{xalignat*}

Furthermore, since the function values and the first $m$ derivatives of $L^{m+1}_{j+m-1}$ vanish in the end points, there holds
\begin{gather*}
     \nodal^0_i(\phi^0_j) = 0, 
     \qquad i=1,\dots,2m+1,
     \quad j=2m+2,\dots,n.
\end{gather*}
What is missing are the node functionals in the second group. If applied to the functions in the first group, there may be nonzero entries anywhere. For the functions in the second group, we have to expand the integrand with respect to Legendre polynomials.
Since there holds
\begin{gather}
    \label{eq:Legendre-integral}
    2(2j+1) L^1_j(x) = \ell_{j+1}(x) - \ell_{j-1}(x),
\end{gather}
we obtain by recursion that
\begin{gather}
    \label{eq:Legendre-expansion}
    L^m_{j} = \alpha_1 \ell_{j+m} + \alpha_2 \ell_{j+m-2} + \cdots + \alpha_{m+1} \ell_{j-m},
\end{gather}
where $\alpha_i$ are coefficients computable from~\eqref{eq:Legendre-integral}.
Using that $L^m_{j+m-1}$ is the derivative of $\phi^0_{2m+j}$ and $\ell_{j-1}$ is the lowest order polynomial in its expansion, we note for $i<j$
\begin{gather*}
 \nodal^0_{2m+i}(\phi^0_{2m+j})= \int_0^1 \ell_{i-1} L^m_{j+m-1} = 0.
\end{gather*}
This yields a lower triangular matrix.
Ignoring the finer structure in formula~\eqref{eq:Legendre-expansion}, we obtain the following block matrix structure:
\begin{gather}
    \label{eq:nodal-matrix}
    \nodal^0_i(\phi^0_j) = \left(
    \begin{array}{ccc|ccc|c}
         * &&&&&&\\
         & \ddots &&&&&\\
         && * &&&& \\\hline
         * & \cdots & * & * &&& \\
         \vdots & & \vdots & \vdots & \ddots & & \\
         * & \cdots & * & * & \cdots & * & \\\hline
         &&&&&&*
    \end{array}
    \right).
\end{gather}
Since this matrix is lower triangular with nonzero diagonal entries, it has full rank and thus the FE $\P_{\deg}\Lambda^0$ is unisolvent. In particular, the set of basis functions and the set of node functionals are both linearly independent.

Furthermore, it has the same block structure as $M_k$ in~\eqref{eq:block-mk}, the kernel being the one dimensional subspace corresponding to the last row and column, while the matrix $A_k$ corresponds to the first $2\times2$-block system.
Due to~\eqref{eq:nodal-matrix-commute}, the matrix $\nodal^1_i(\phi^1_j)$ is obtained by deleting the last row and column from~\eqref{eq:nodal-matrix}. Therefore, $\P_{\deg-1}\Lambda^1$ is unisolvent as well. 

%%%%%%%%%%%%%%%%%%%%%%%%%%%%%%%%%%%%%%%%%%%%%%%%%%%%%%%%%%%%%%%%%%%%%%%%%%%%%%%%%%%%%%%%
\subsection{Commutation property of interpolation operators}
%%%%%%%%%%%%%%%%%%%%%%%%%%%%%%%%%%%%%%%%%%%%%%%%%%%%%%%%%%%%%%%%%%%%%%%%%%%%%%%%%%%%%%%%

The node functionals and shape functions in subsections~\ref{sec:fem0} and~\ref{sec:fem1} serve to define interpolation operators according to~\eqref{eq:cint}.
They are constructed precisely to fulfil the assumptions of Lemma~\ref{lem:commuting-lemma}.
Thus, by the use of this construction, we obtain a commuting interpolation operator.

\begin{remark}
While the node functionals in~\eqref{eq:nodes-1} are standard degrees of freedom for Hermite interpolation, those in~\eqref{eq:nodes-2} and~\eqref{eq:nodes-3} are linear combinations of the usual ones. More so, they make the implementation considerably more complicated, since their usual purpose of ensuring continuity is not easily accomplished. Therefore, we point out that these degrees of freedom serve for the analysis of interpolation operators only, but that implementations should use $u(0)$ and $u(1)$ as in Hermite interpolation.
\end{remark}

\begin{remark}
The node functionals $\{\nodal^0_i\}_{i=1,\ldots,n+1}$ and $\{\nodal^1_i\}_{i=1,\ldots,n}$ are well-defined for $C^m$-reg\-u\-lar functions. Following the same procedure as in~\cite{BonizzoniKanschat21}, it is possible to weaken such conditions, and introduce weighted node functionals $\{\overline{\nodal^0_i}\}_{i=1,\ldots,n+1}$ and $\{\overline{\nodal^1_i}\}_{i=1,\ldots,n}$ which are bounded in $L^2$. This construction yields the $L^2$-stable and commuting quasi-interpolation operators given by:
\begin{gather}
    \label{eq:quasi_interp}
     \wint_\nu(u) = \sum_{i=1}^{n+1} \overline{\nodal^\nu_i}(u) \left(\sum_{j=1}^{n_\nu} \alpha^\nu_{ij}\phi^\nu_j\right),
\end{gather}
for $\nu=0, 1$.
\end{remark}

%%%%%%%%%%%%%%%%%%%%%%%%%%%%%%%%%%%%%%%%%%%%%%%%%%%%%%%%%%%%
\section{Tensorization of the Finite Element Complex}
%%%%%%%%%%%%%%%%%%%%%%%%%%%%%%%%%%%%%%%%%%%%%%%%%%%%%%%%%%%%

In the present section we apply the tensor product construction of cochain complexes (see~\cite{ArnoldBoffiBonizzoni}) to the FE de Rham subcomplex~\eqref{eq:1D_complex}. Moreover, tensor product interpolation operators are introduced.

%%%%%%%%%%%%%%%%%%%%%%%%%%%%%%%%%%%%%%%%%%%%%%%%%%%%%%%%%%%%
\subsection{Tensor complex in two dimensions}
%%%%%%%%%%%%%%%%%%%%%%%%%%%%%%%%%%%%%%%%%%%%%%%%%%%%%%%%%%%%

The tensor product of complex~\eqref{eq:1D_complex} on $\interval$ with itself is the following de Rham complex on the Cartesian product $\interval\times\interval$:
\begin{gather}
\label{eq:tp-complex-2}
    0 \xrightarrow{\;\;\subset \;\;}
    (\P\Lambda^{\otimes 2})^0
    \xrightarrow{\;\; \dif \;\;}
    (\P\Lambda^{\otimes 2})^1
    \xrightarrow{\;\; \dif \;\;}
    (\P\Lambda^{\otimes 2})^2
    \xrightarrow{\;\;0 \;\;}0
\end{gather}
where the spaces of polynomial differential forms are defined as
\begin{equation*}
\begin{aligned}
    (\P\Lambda^{\otimes 2})^0 & = \P_n\Lambda^0 \otimes \P_n\Lambda^0,\\
    (\P\Lambda^{\otimes 2})^1 & = \left(\P_n\Lambda^0 \otimes \P_{n-1}\Lambda^1\right)\oplus\left(\P_{n-1}\Lambda^1 \otimes \P_n\Lambda^0\right),\\
    (\P\Lambda^{\otimes 2})^2 & = \P_{n-1}\Lambda^1 \otimes \P_{n-1}\Lambda^1,
\end{aligned}
\end{equation*}
and the exterior derivative $\dif\colon(\P\Lambda^{\otimes 2})^\nu\rightarrow(\P\Lambda^{\otimes 2})^{\nu+1}$, $\nu=0,1,2$, is the linear operator acting on tensor product functions as follows:
\begin{equation}
\label{eq:d_tens2}
\begin{aligned}
    \dif (u_0\otimes v_0) &= \dif u_0\otimes v_0 + u_0\otimes \dif v_0,\quad
    &\forall\, u_0&\in \P_n\Lambda^0,
    & v_0&\in \P_n\Lambda^0,\\
    \dif (u_0\otimes v_1) &= \dif u_0\otimes v_1,
    &\forall\, u_0&\in \P_n\Lambda^0,
    & v_1&\in \P_{n-1}\Lambda^1,\\
    \dif (u_1\otimes v_0) &= u_1\otimes \dif v_0,
    & \forall\, u_1&\in \P_{n-1}\Lambda^1,
    & v_0&\in \P_n\Lambda^0,\\
    \dif (u_1\otimes v_1) &= 0,
    & \forall\, u_1&\in \P_{n-1}\Lambda^1,
    & v_0&\in \P_{n-1}\Lambda^1.
\end{aligned}
\end{equation}
We note that~\eqref{eq:tp-complex-2} is a subcomplex of the de Rham complex of smooth forms on $\interval\times\interval$, as $(\P\Lambda^{\otimes 2})^\nu\subset C^\infty\Lambda^\nu(\interval\times\interval)$.

A basis for $(\P\Lambda^{\otimes 2})^\nu$, $\nu=0,1,2$, can be obtained by tensorization of the basis $\{\varphi_j^0\}_{j=1}^{n+1}$ and $\{\varphi^1_j\}_{j=1}^n$ of $\P_n\Lambda^0$ and $\P_{n-1}\Lambda^1$, respectively, namely
\begin{equation}
\begin{aligned}
    \label{eq:tp-basis-2}
    (\P\Lambda^{\otimes 2})^0 & 
    = {\rm span}\left\{\varphi^0_{j_1}\otimes\varphi^0_{j_2},\, 1\leq j_1, j_2\leq n+1\right\},\\
    (\P\Lambda^{\otimes 2})^1 & 
    = {\rm span}\left\{\varphi^0_{j_1}\otimes\varphi^1_{j_2},\, 1\leq j_1\leq n+1,\; 1\leq j_2\leq n\right\} \\
    &\oplus
    {\rm span}\left\{\varphi^1_{j_1}\otimes\varphi^0_{j_2},\, 1\leq j_1\leq n,\; 1\leq j_2\leq n+1\right\} ,\\
    (\P\Lambda^{\otimes 2})^2 & 
    = {\rm span}\left\{\varphi^1_{j_1}\otimes\varphi^1_{j_2},\, 1\leq j_1, j_2\leq n\right\}.
\end{aligned}
\end{equation}

The node functionals for $(\P\Lambda^{\otimes 2})^\nu$, $\nu=0,1,2$, are defined by tensor product and apply to tensor product functions as follows:
\begin{equation}
    \label{eq:node-functionals-2d}
    \nodal^{(i_1,i_2)}_{(j_1,j_2)}
        (u\otimes v)
    =\left[\nodal^{i_1}_{j_1}\otimes\nodal^{i_2}_{j_2}\right]
        (u\otimes v)
    = \nodal^{i_1}_{j_1}(u)
    \nodal^{i_2}_{j_2}(v),
\end{equation}
for all compatible pairs of indices $(i_1,i_2)$ and $(j_1,j_2)$. For example, in the particular case $n=5$, the $C^2$-regular tensor product FEs are complemented with the following node functionals:
\begin{itemize}
    \item for $k=0$ and $u,\, v\in\P_5\Lambda^0$
    \begin{xalignat*}2
        \nodal^{00}_{11}(u\otimes v)&=\partial_x u(0)\partial_y v(0)&
        \nodal^{00}_{12}(u\otimes v)&=\partial_x u(0)\partial_y v(1)\\
        \nodal^{00}_{13}(u\otimes v)&=\partial_x u(0)\partial^2_y v(0)&
        \nodal^{00}_{14}(u\otimes v)&=\partial_x u(0)\partial^2_y v(1)\\
        \nodal^{00}_{15}(u\otimes v)&=\partial_x u(0)(v(1)-v(0))&
        \nodal^{00}_{16}(u\otimes v)&=\partial_x u(0)(v(1)+v(0))
    \end{xalignat*}
    \item for $k=1$ and $u\in\P_5\Lambda^0$, $v\in\P_{4}\Lambda^1$
    \begin{xalignat*}2
        \nodal^{01}_{11}(u\otimes v)&=\partial_x u(0)v(0)&
        \nodal^{01}_{12}(u\otimes v)&=\partial_x u(0)v(1)\\
        \nodal^{01}_{13}(u\otimes v)&=\partial_x u(0)\partial_y v(0)&
        \nodal^{01}_{14}(u\otimes v)&=\partial_x u(0)\partial_y v(1)\\
        \nodal^{01}_{15}(u\otimes v)&=\partial_x u(0)\int_0^1 v\,\dx&&
    \end{xalignat*}
    \item for $k=2$ and $u,\, v\in\P_{4}\Lambda^1$
    \begin{xalignat*}2
        \nodal^{11}_{11}(u\otimes v)&=u(0)v(0)&
        \nodal^{11}_{12}(u\otimes v)&=u(0)v(1)\\
        \nodal^{11}_{13}(u\otimes v)&=u(0)\partial_y v(0)&
        \nodal^{11}_{14}(u\otimes v)&=u(0)\partial_y v(1)\\
        \nodal^{11}_{15}(u\otimes v)&=u(0)\int_0^1 v\, \dx&&
    \end{xalignat*}
\end{itemize}
Tensor product node functionals apply to two-dimensional functions in a straightforward way. For example, for $u\in\mathcal C^\infty\Lambda^0(\interval\times\interval)$, there holds:
\begin{xalignat*}2
        \nodal^{00}_{11}(u)&=\partial_x\partial_y u(0,0)&
        \nodal^{00}_{12}(u)&=\partial_x\partial_y u(0,1)\\
        \nodal^{00}_{13}(u)&=\partial_x\partial^2_y u(0,0)&
        \nodal^{00}_{14}(u)&=\partial_x\partial^2_y u(0,1)\\
        \nodal^{00}_{15}(u)&=\partial_x u(0,1) - \partial_x u(0,0)&
        \nodal^{00}_{16}(u)&=\partial_x u(0,1) + \partial_x u(0,0)
\end{xalignat*}

Finally, the tensor product interpolation operator $\cint^{\otimes 2}_\nu$, $\nu=0,1,2$, is defined as follows:
\begin{equation*}
\begin{aligned}
    \cint^{\otimes 2}_0&=\cint_0\otimes \cint_0,\\
    \cint^{\otimes 2}_1&=\cint_0\otimes \cint_1 + \cint_1\otimes \cint_0,\\
    \cint^{\otimes 2}_2&=\cint_1\otimes \cint_1,
\end{aligned}
\end{equation*}
It applies to rank-one functions of the form $u_1\otimes u_2\in \Lambda^{i_1}\otimes\Lambda^{i_2}$ as
\begin{equation*}
    \cint^{\otimes 2}_\nu(u_1\otimes u_2)
    =I_{i_1}(u_1) \otimes I_{i_2}(u_2),
\end{equation*}
and can be naturally applied to two-dimensional functions.

%%%%%%%%%%%%%%%%%%%%%%%%%%%%%%%%%%%%%%%%%%%%%%%%%%%%%%%%%%%%
\subsection{Tensor complex in $\ordertens$ dimensions}
%%%%%%%%%%%%%%%%%%%%%%%%%%%%%%%%%%%%%%%%%%%%%%%%%%%%%%%%%%%%

The tensor product construction generalizes to any tensorization order $\ordertens\geq 2$, giving rise to the following de Rham subcomplex on the Cartesian product $\interval^{\times\ordertens}$:
\begin{gather*}
    0 \xrightarrow{\;\;\subset \;\;}
    \cdots
    \xrightarrow{\;\; \dif \;\;}
    (\P\Lambda^{\otimes \ordertens})^\nu
    \xrightarrow{\;\; \dif \;\;}
    (\P\Lambda^{\otimes \ordertens})^{\nu+1}
    \xrightarrow{\;\; \dif \;\;}
    \cdots
    \xrightarrow{\;\;0 \;\;}0.
\end{gather*}
The tensor product spaces are defined as
\begin{gather}
    \label{eq:Ptens}
    (\P\Lambda^{\otimes \ordertens})^\nu
    = \bigoplus_{\mathbf i\in\chi_\nu}
    \left(\P_{n-i_1}\Lambda^{i_1}\otimes\cdots\otimes\P_{n-i_\ordertens}\Lambda^{i_\ordertens}\right)\subset C^\infty\Lambda^\nu(\interval^{\times N}),
\end{gather}
where $\interval^{\times N}=\underbrace{\interval\times\cdots\times\interval}_{N\rm{ times}}$ and $\chi_\nu$ is the characteristic vector of a combination $\nu$ out of $N$ given by 
\begin{gather}
\label{eq:chi}
    \chi_\nu:=\left\{\boldsymbol{i}=(i_1,\ldots,i_\ordertens)\in\{0,1\}^\ordertens
    \;\middle|\;\sum_{\ell=1}^\ordertens i_\ell=\nu\right\}.
\end{gather}
In analogy to formula~\eqref{eq:tp-basis-2}, there holds
\begin{equation}
\label{eq:tp-basis-n}
    (\P\Lambda^{\otimes \ordertens})^\nu
    =\bigoplus_{\mathbf i\in\chi_\nu}{\rm span}\{\varphi^{\mathbf i}_{\mathbf j},\, 1\leq j_\ell\leq n+1-i_\ell,\, \ell=1,\ldots,\ordertens
    \},
\end{equation}
with $\varphi^{\mathbf i}_{\mathbf j}\colon=\varphi^{i_1}_{j_1}\otimes\cdots\otimes\varphi^{i_\ordertens}_{j_\ordertens}$. The tensor product basis $\{\varphi^{\mathbf i}_{\mathbf j},\, \mathbf{i}\in\chi_\nu,\, \mathbf{j}\leq n+1-\mathbf{i}\}$ is also known as rank-one basis.

The exterior derivative acts on the rank-one basis function $\varphi^{\mathbf{i}}_{\mathbf{j}}\in(\P\Lambda^{\otimes \ordertens})^\nu$ as
\begin{equation}
\label{eq:d-tensor-k}
    \dif \varphi^{\mathbf{i}}_{\mathbf{j}}
    = \dif (\varphi^{i_1}_{j_1}\otimes \cdots\otimes \varphi^{i_\ordertens}_{j_\ordertens})
    = \sum_{t=1}^\ordertens \theta_t\,
    (\varphi^{i_1}_{j_1}\otimes\cdots\otimes 
    \dif \varphi^{i_t}_{j_t}\otimes\cdots\otimes \varphi^{i_\ordertens}_{j_\ordertens}),
\end{equation}
with $\theta_t:=(-1)^{\sum_{\ell=1}^{j-1} i_\ell}\in\{-1,1\}$.
Formula~\eqref{eq:d-tensor-k} extends by linearity to any element of $(\P\Lambda^{\otimes \ordertens})^\nu$.

% %%%%%%%%%%%%%%%%%%%%%%%%%%%%%%%%%%%%%%%%%%%%%%%%%%%%%%%%%%%%
% \subsection{Tensor product of interpolation operators}
% %%%%%%%%%%%%%%%%%%%%%%%%%%%%%%%%%%%%%%%%%%%%%%%%%%%%%%%%%%%%

The tensor product construction yields node functionals $\{\nodal^{\mathbf i}_{\mathbf j},\, \mathbf{i}\in\chi_\nu,\,1\leq j\leq n+1-\mathbf i\}$ for $(\P\Lambda^{\otimes \ordertens})^\nu$ that are defined on rank-one basis functions as
\begin{equation}
    \label{eq:node-functionals-tp}
    \nodal^{\mathbf i}_{\mathbf j}(\varphi^{\mathbf{i}}_{\mathbf{j}^\prime})
    =\nodal^{i_1}_{j_1}\otimes\cdots\otimes\nodal^{i_\ordertens}_{j_\ordertens}
    (\varphi^{i_1}_{j_1^\prime}\otimes \cdots\otimes \varphi^{i_\ordertens}_{j_\ordertens^\prime})
    = \nodal^{i_1}_{j_1}(\varphi^{i_1}_{j_1^\prime})\cdots
    \nodal^{i_\ordertens}_{j_\ordertens}(\varphi^{i_\ordertens}_{j_\ordertens^\prime}),
\end{equation}
and then extend by linearity to any element of $(\P\Lambda^{\otimes \ordertens})^\nu$, with the convention $\nodal^{\mathbf i}_{\mathbf j}(\varphi^{\mathbf{i}^\prime}_{\mathbf{j}^\prime})=0$ for $\mathbf i\neq\mathbf{i}^\prime$.
The node functionals apply to higher-dimensional regular functions in a straightforward way. 

Having introduced node functionals and rank-one basis functions for $(\P\Lambda^{\otimes \ordertens})^\nu$, it is natural to define the $\ordertens$-dimensional interpolation operator, which coincides with the tensor product interpolation operator  
\begin{equation}
    \label{eq:tp-interp}
    \cint^{\otimes \ordertens}_\nu:=\sum_{\mathbf i\in\chi_\nu} \cint_{i_1}\otimes\cdots\otimes\cint_{i_\ordertens},
\end{equation}
once we extend $\cint_0$ on 1-forms and $\cint_1$ on 0-forms as the null-operator.

Exploiting the tensor product construction, and using Lemma~\ref{lem:commuting-lemma}, the following result can be proved.  
\begin{lemma}
\label{lem:tp-interp-cochain}
The tensor product interpolation operator $\cint^{\otimes\ordertens}_\nu$ defined in~\eqref{eq:tp-interp} is a cochain operator, namely, for all smooth forms $u\in C^\infty\Lambda^\nu(\interval^{\times\ordertens})$, there holds
  \begin{equation}
    \label{eq:tp-interp-cochain}
      \cint^{\otimes\ordertens}_{\nu+1} \dif u
      = \dif \cint^{\otimes\ordertens}_{\nu} u.
  \end{equation}
\end{lemma}
For a complete proof we refer the reader to~\cite{BonizzoniKanschat21}, where equation~\eqref{eq:tp-interp-cochain} is first proven for rank-one functions, and then extended by linearity and density to all smooth differential forms. 

\begin{remark}
The tensor product construction applied to the quasi-interpolation operators~\eqref{eq:quasi_interp} yields tensor product quasi-interpolation operators that are cochain operators and bounded in $L^2$. A similar tensor product construction has been applied to projection operators in~\cite{BonizzoniBuffaNobile13}.
\end{remark}

%%%%%%%%%%%%%%%%%%%%%%%%%%%%%%%%%%%%%%%%%%%%%%%%%%%%%%%%%%%%
\section{Conclusions and Further Developments}
%%%%%%%%%%%%%%%%%%%%%%%%%%%%%%%%%%%%%%%%%%%%%%%%%%%%%%%%%%%%

We have developed tensor product FEs for the de Rham complex of arbitrary smoothness and dimension. Their construction is based on one dimensional Hermite interpolation and the commutation property of the finite element interpolation operators is proven using a general commutation lemma.

This construction plays an important role in the framework of Bernstein-Gelfand-Gelfand (BGG) sequences~\cite{ArnoldHu}, where for instance de Rham complexes of different smoothness are combined to yield new complexes. The construction here yields the one-dimensional building blocks for tensor product FE BGG sequences as well as the tensor product base complexes, as we point out in~\cite{BonizzoniHuKanschatSap22}.

The general framework of de Rham complex embeds several problems important from the application point of view, like the Darcy and the Maxwell problem. To discretize those problems, it is then possible to employ the proposed FE spaces. In particular, they can be combined with model order reduction techniques, to handle parametric boundary value problems. Some works in this direction~\cite{BonizzoniNobilePerugia18,BonizzoniNobilePerugiaPradovera20a,BonizzoniNobilePerugiaPradovera20b,BonizzoniPradovera,BonizzoniPradoveraRuggeri} consider the parametric Helmholtz equation, and they could represent the starting point towards the parametric Maxwell equation. Moreover, following the ideas in~\cite{BonizzoniBuffaNobile13}, they can be employed in combination with uncertainty quantification techniques to deal with the lognormal Darcy problem modeling the fluid flow in bounded heterogeneous media~\cite{BonizzoniNobile20, BonizzoniNobileKressner,BonizzoniNobile14,BonizzoniNobile12}.

\section*{Acknowledgements}
The first author acknowledge support from the HGS MathComp through the Distinguished Romberg Guest Professorship Program.
Moreover, the work of the first author is part of a project that has received funding from the European Research Council ERC under the European Union's Horizon 2020 research and innovation program (Grant agreement No.~865751). The second author was supported by the Deutsche Forschungsgemeinschaft (DFG, German Research Foundation) under Germany's Excellence Strategy EXC 2181/1 - 390900948 (the Heidelberg STRUCTURES Excellence Cluster).


\begin{thebibliography}{99}

\bibitem{ArnoldBoffiBonizzoni}
Arnold, D. N. and Boffi, D. and Bonizzoni, F.
Finite element differential forms on curvilinear cubic meshes and their approximation properties.
\textit{Numerische Mathematik} (2015) \textbf{129} (1): 1--20.

\bibitem{ArnoldFalkWinther}
Arnold, D. N. and Falk, R. S. and Winther, R.
Finite element exterior calculus, homological techniques, and
applications.
\textit{Acta Numerica} (2006) 15: 1--155.

\bibitem{ArnoldHu}
Arnold, D. N. and Hu, K.
Complexes from Complexes
\textit{Found. Comput. Math.} (2021) \textbf{21} (6): 1739--1774.

\bibitem{BonizzoniBuffaNobile13}
Bonizzoni, F. and Buffa, A. and Nobile, F.
Moment equations for the mixed formulation of the Hodge Laplacian with stochastic loading term.
\textit{IMA Journal of Numerical Analysis} (2013) \textbf{34} (4): 1328--1360.

\bibitem{BonizzoniKanschat21}
Bonizzoni, F. and Kanschat, G.
{$H^1$}-conforming finite element cochain complexes and commuting quasi-interpolation operators on Cartesian meshes.
\textit{Calcolo} (2021) \textbf{58} (2): 1--29.

\bibitem{BonizzoniHuKanschatSap22}
Bonizzoni, F. and Hu, K. and Kanschat, G. and Sap, D.
{Spline and tensor product finite element BGG sequences}
\textit{In preparation} (2022).

\bibitem{BonizzoniNobile12}
Bonizzoni, F. and Nobile, F.
Perturbation analysis for the stochastic Darcy problem.
\textit{Proceeding in ECCOMAS 2012-European Congress on Computational Methods in Applied Sciences and Engineering} (2012): 3926--3933.
ISBN: 9783950353709

\bibitem{BonizzoniNobile14}
Bonizzoni, F. and Nobile, F.
Perturbation analysis for the Darcy problem with Log-normal permeability.
\textit{SIAM/ASA Journal on Uncertainty Quantification} (2014) \textbf{2} (1): 223 -- 244.

\bibitem{BonizzoniNobile20}
Bonizzoni, F. and Nobile, F.
Regularity and sparse approximation of the recursive first moment equations for the lognormal Darcy problem.
\textit{Computers and Mathematics with Applications} (2020) \textbf{80} (12): 2925 -- 2947.

\bibitem{BonizzoniNobileKressner}
Bonizzoni, F. and Nobile, F. and Kressner, D.
Tensor Train approximation of moment equations for elliptic equations with lognormal coefficient.
\textit{Computer Methods in Applied Mechanics and Engineering} (2016) \textbf{308}: 349 -- 376.

\bibitem{BonizzoniNobilePerugia18}
Bonizzoni, F. and Nobile, F. and Perugia, I.
Convergence analysis of Pad\'e approximations for Helmholtz frequency response problems.
\textit{ESAIM: Mathematical Modelling and Numerical Analysis} (2018)  \textbf{52} (4): 1261 -- 1284.

\bibitem{BonizzoniNobilePerugiaPradovera20a}
Bonizzoni, F. and Nobile, F. and Perugia, I. and Pradovera, D.
Fast Least-Squares Pad\'e approximation of problems with normal operators and meromorphic structure.
\textit{Mathematics of Computation} (2020) \textbf{89}: 1229--1257.

\bibitem{BonizzoniNobilePerugiaPradovera20b}
Bonizzoni, F. and Nobile, F. and Perugia, I. and Pradovera, D.
Least-Squares Pad\'e approximation of parametric and stochastic Helmholtz maps.
\textit{Advances in Computational Mathematics} (2020) \textbf{46} (46).

\bibitem{BonizzoniPradovera}
Bonizzoni, F. and Pradovera, D.
Shape optimization for a noise reduction problem by non-intrusive parametric reduced modeling. 
\textit{Proceeding in the 14th WCCM-ECCOMAS Congress 2020/2021}, (2021).
DOI: 10.23967/wccm-eccomas.2020.300

\bibitem{BonizzoniPradoveraRuggeri}
Bonizzoni, F. and Pradovera, D. and Ruggeri, M.
Rational-based model order reduction of Helmholtz frequency response problems with adaptive finite element snapshots.
(2022) \textit{arXiv:2112.04302}.

\end{thebibliography}
\end{document}